\documentclass[12pt]{article}
\usepackage{amssymb}
\usepackage{amsthm}
\usepackage{amsmath}
\usepackage{graphicx}
\graphicspath{{./Figs/}}
\usepackage{lmodern}
\usepackage{graphics}
\usepackage{float}

\newtheorem{teo}{Theorem}
\newtheorem{lem}{Lemma}

\newtheorem{remark}{Remark}

\newcommand{\pa}{\partial}

\newcommand{\ve}{\varepsilon}
\newcommand{\om}{\omega}
\newcommand{\be}{\begin{equation}}
\newcommand{\ee}{\end{equation}}

\newcommand{\ipd}{\stackrel{\normalfont\text{def}}{=}}
\newcommand{\const}{\operatorname{const}}

\usepackage[usenames]{color}
\begin{document}
\allowdisplaybreaks

\title{Solitons and peaked solitons for the general Degasperis-Procesi model}
\author{J.~Noyola Rodriguez
\thanks{Universidad de Sonora, Mexico,\
jesnoyola89@gmail.com}\and  G.~Omel'yanov\thanks{
Universidad de Sonora, Mexico,\ omel@mat.uson.mx} }

\date{}
\maketitle
\begin{abstract}
We consider  the general Degasperis-Procesi model of shallow water out-flows. This six parametric family of conservation laws contains, in particular,   KdV,  Benjamin-Bona-Mahony, Camassa-Holm, and Degasperis-Procesi    equations.  The main result consists of  criterions which guarantee  the existence of solitary wave solutions: solitons and peakons ("peaked solitons").
\end{abstract}



\section{Introduction}
 We consider a modern unidirectional approximation of the shallow water system called
  the  ``general Degasperis-Procesi" model (\cite{ DegProc}, 1999):
\begin{align}
&\frac{\pa }{\pa t}\left\{u-\alpha^2\ve^2\frac{\pa^2 u}{\pa x^2}\right\}\label{1}\\
&+\frac{\pa}{\pa x}\left\{c_0u+c_1u^2
-c_2\ve^2\Big(\frac{\pa u}{\pa x}\Big)^2+\ve^2\big(\gamma-c_3u\big)\frac{\pa^2 u}{\pa x^2}\right\}=0, \; x \in \mathbb{R}^1, \; t > 0.\notag
\end{align}
Here $\alpha$, $c_0,\dots,c_3$, $\gamma$ are  real parameters  and  $\varepsilon$ characterizes the dispersion.

This six parametric family of third order conservation laws contains as particular cases a list of basic equations. Indeed:

1. If we set $\alpha=c_2=c_3=0$ then we obtain the famous KdV equation, whereas for
 $\gamma=c_2=c_3=0$ Eq.(\ref{1}) is the well known  Benjamin-Bona-Mahony (BBM) equation (\cite{BBM}, 1972).

2. Preserving in  (\ref{1}) the nonlinear dispersion terms     and setting
  $c_2=c_3/2$, $c_1=3c_3/2\alpha^2$, and $\gamma=0$ we obtain the Camassa-Holm (CH) equation (\cite{CH}, 1993):
\begin{equation}
\frac{\pa }{\pa t}\left\{u-\alpha^2\ve^2\frac{\pa^2 u}{\pa x^2}\right\}
+\frac{\pa}{\pa x}\left\{c_0u+c_1u^2
-\ve^2\frac{c_3}{2}\Big\{\Big(\frac{\pa u}{\pa x}\Big)^2+2u\frac{\pa^2 u}{\pa x^2}\Big\}\right\}=0. \label{3}
\end{equation}
3. In the case $c_2=c_3$, $c_1=2c_3/\alpha^2$, and $c_0=\gamma=0$ (\ref{1}) is the Degasperis-Procesi (DP) equation (\cite{DegProc}, see also \cite{ELY} and references therein):
\begin{equation}
\frac{\pa }{\pa t}\left\{u-\alpha^2\ve^2\frac{\pa^2 u}{\pa x^2}\right\}
+\frac{\pa}{\pa x}\left\{c_1u^2
-\ve^2c_3\Big\{\Big(\frac{\pa u}{\pa x}\Big)^2+u\frac{\pa^2 u}{\pa x^2}\Big\}\right\}=0. \label{4}
\end{equation}

The KdV and BBM equations are essentially different. Both of them have soliton-type traveling wave solutions, however, KdV solitons collide elastically: they pass through each other preserving the shape and velocities, whereas BBM ``solitons" have changed after the interaction and an oscillatory tail is generated \cite{BPS}.

Next,
for the first view the CH (\ref{3}) and DP (\ref{4}) equations are quite similar: the difference consists of the relation between the coefficients $c_2$ and $c_3$ only.
However, it should be emphasized that these equations have truly different properties:

- if $c_0>0$, the Camassa-Holm equation has smooth soliton solution
\begin{equation}\label{02}
u(x,t)=A\om\big((x-Vt)/\varepsilon, A\big),\quad \om(\tau,\cdot)\in C^\infty(\mathbb{R}^1),\quad \lim_{\tau\to\pm\infty}\om(\tau,\cdot)=0,
\end{equation}
- if $c_0=0$,  the Camassa-Holm equation has the so-called "peakon" solitary wave solution \cite{CH}, that is a continuous function of the form
\be
u(x,t)= A\exp\big(-|x-Vt|/\ve\big),\quad V>0,\label{03}
\ee
- the Degasperis-Procesi equation, under the condition $u\to0$ as $x\to\pm\infty$, has non-smooth traveling wave solutions only. Namely, peakon-type solution of the form (\ref{03}) and "shock-peakon" \cite{DHH}, which is given by
\be
u(x,t)= -(t+k)sgn(x-Vt)\exp\big(-|x-Vt|/\ve\big),\quad V>0,\quad k>0.\label{04}
\ee
Note also that there are many other solutions of (\ref{4}) if we alow $u\to\const\neq0$ as $x\to\pm\infty$ \cite{ZQQ}.

To justify peakon as a well-defined solution of DP equation  we transform (\ref{1}) to the following divergent form:
\begin{align}
\frac{\pa u}{\pa t}+ \frac{\pa}{\pa x}&\left\{c_0u+c_1u^2-(c_2-c_3)\Big(\ve\frac{\pa u}{\pa x}\Big)^2 \right\}\label{05}\\
&=\ve^2\frac{\pa^2}{\pa x^2}\left\{\alpha^2\frac{\pa u}{\pa t}-\frac{\pa }{\pa x}\left(\gamma u-\frac{c_3}{2}u^2\right)\right\}, \notag
\end{align}
and note that all terms here are well defined not for smooth functions only, but  for distributions of the type (\ref{03}) also.
As for shock-peakon (\ref{04}), it seems that the Degasperis-Procesi equation (\ref{4}) is the unique representative from the family   (\ref{1}), for which  such type of solutions can be defined correctly. Indeed, (\ref{04}) is the jump-type function, thus $u_x\sim \delta(x-Vt)$ for any $\ve=\const$. Therefore, the term $(u_x)^2$ doesn't exist in the weak sense  and the equation (\ref{05}) with $c_2\neq c_3$ is bad defined in the sense of distributions.

The difference between the equations (\ref{3}) and (\ref{4}) can be demonstrated also by use the balance law for the basic model (\ref{1}):
\begin{equation}
\frac{d }{d t}\left\{\int_{-\infty}^\infty u^2dx+\alpha^2\int_{-\infty}^\infty(\ve u_x)^2dx\right\}=\ve^{-1}(c_3-2c_2)\int_{-\infty}^\infty(\ve u_x)^3dx.\label{06}
\end{equation}
It is clear that the Camassa-Holm equation with $c_3=2c_2$ is the exclusive situation when (\ref{06}) is the conservation law, whereas all other relations between $c_2$ and $c_3$ imply, generally speaking, instability of the solution. The Cauchy problems for the CH and DP equations have been studied extensively. We  refer  readers to the paper  \cite{ELY}, which contains further references also.

Three particular cases, i.e. the equations KdV, CH (\ref{3}), and DP (\ref{4}) belong to the so-called "integrable equations" (see e.g. \cite{DegProc, CH, DHH}, \cite{Con}-\cite{Mat}). In particular, it is known that the solitary waves  interact elastically in these models. At the same time, returning to the gDP model, we stress that these special cases exhaust that's all what is known about the general family (\ref{1}). In particular, it remains unknown how to divide the space of structural parameters $\alpha$, $\gamma$, $c_i$ in order to separate smooth and non-smooth traveling wave solutions. Furthermore, excepting the KdV, CH, and DP equations; all other versions of the model (\ref{1}) are essentially non-integrable (see e.g. \cite{ELY}). Respectively, the character of wave collision remains unknown also.

 To begin the study of wave propagations for non-integrable versions of (\ref{1})  we should separate firstly two basic situations: smooth and non-smooth traveling solutions. Section 2 contains the construction of solitons and obtaining sufficient conditions for their existence.  The non-smooth case is considered in Section 3.  We use an alternative approach there and show that peakons are just peaked solitons (see also \cite{Mat}).
\section{Soliton type solution}
Similar to (\ref{02}) let us set the ansatz
\begin{equation}\label{5}
u=A\om\big(\beta(x-Vt)/\varepsilon,A\big),
\end{equation}
where $\om(\eta,A)$ is a smooth  function such that
\begin{align}
&\om(-\eta,A)=\om(\eta,A),\quad\om(\eta,A)\to0\quad\text{as}\quad \eta\to\pm\infty,\label{6}\\
&\om(0,A)=1,\label{7}
\end{align}
the amplitude $A>0$ is a  free parameter, and the velocity
 $V=V(A)$ should be determined. To simplify formulas we define  the scale $\beta=\sqrt{c_1(c_2+c_3)}/c_3$.

 In what follows we assume that
\be\label{8}
\gamma\geq0,\quad \alpha\geq0, \quad \gamma+\alpha>0,\quad c_0\geq0,\quad c_k>0,\quad k=1,2,3.
\ee
Substituting  (\ref{5}) into Eq.(\ref{1}), integrating, and using (\ref{6}), we obtain the following version of the inverse scattering problem:

\textit{Determine the velocity} $V$ \textit{such that the equation}
\begin{align}
\Big\{1-&\frac{c_3A }{\gamma+\alpha^2V}\om\Big\}\frac{d^2 \om}{d \eta^2}=\frac{c_2A }{\gamma+\alpha^2V}\left(\frac{d \om}{d \eta}\right)^2\notag\\
&+\frac{c_3^2}{c_1(c_2+c_3)(\gamma+\alpha^2V)}\big((V-c_0)\om-c_1A\om^2\big),\quad \eta\in \mathbb{R}^1,\label{9}
\end{align}
\textit{has a nontrivial smooth solution with the properties} (\ref{6}), (\ref{7}).

Let us simplify this problem. To this end rescaling the function $\om$,
\begin{equation}\label{11}
 W=p\om,\quad p=c_3A/(\gamma+\alpha^2V),
\end{equation}
we  define
\begin{equation}
r=c_3/(c_2+c_3), \quad q=c_3(V-c_0)/\big(c_1(\gamma+\alpha^2V)\big),\label{12}
\end{equation}
and pass to the equation
\begin{equation}\label{13}
(1-W)\frac{d^2 W}{d\eta^2}=\frac{1-r}{r}\left(\frac{d W}{d\eta}\right)^2+r(q\,W-W^2),\quad \eta\in \mathbb{R}^1.
\end{equation}
The next step is the substitution
\begin{equation}\label{14}
 W(\eta)=1-g(\eta)^r,
\end{equation}
which allows us to eliminate the first derivatives from the model equation (\ref{13}).
We take into account the condition (\ref{7}) and the property of being even, $g(-\eta)=g(\eta)$. Then under the condition
\be
p<1\label{17}
\ee
 we can pass from the inverse scattering problem (\ref{9}) to the "boundary" problem
\begin{align}
&\frac{d^2 g}{d\eta^2}=g-(2-q)g^{1-r}+(1-q)g^{1-2r},\quad \eta\in(0,\infty),\label{150}\\
& g^r\big|_{\eta=0}=1-p,\quad g|_{\eta\to\infty}=1,\quad dg/d\eta\big|_{\eta=0}=0.\label{16}
\end{align}
 Now we integrate (\ref{150}) and obtain  the first order ODE
\be
\frac{d g}{d\eta}=\sqrt{F(g,q)},\quad \eta\in(0,\infty);\quad
 g|_{\eta=0}=g_*,\label{18}
 \ee
 supplemented by the condition
 \be
F(g_*,q)=0.\label{118}
\ee
Here
\begin{align}
&F(g,q)=g^2-2\frac{2-q}{2-r}g^{2-r}+\frac{1-q}{1-r}g^{2-2r}-C(q),\label{19}\\
&C(q)=r\frac{r-q}{(1-r)(2-r)},\quad g_*=(1-p)^{1/r}.\label{20}
\end{align}
Obviously, the solution of the problem (\ref{18}) with each $q=\const$ exists for $\eta\geq0$, however,  it is unique for $\eta\geq\const>0$ only.

Note now that for each constant $q$
$$
F(1,q)=0,\quad F'_g(1,q)=0, \quad F''_{gg}(1,q)>0.
$$
 Moreover, $F(g,\cdot)$ has only two critical points: $g=1$ and $g=(1-q)^r$. Thus, for all $q$, the condition   (\ref{118}) can be verified not more than at one point $g_*<1$.

Note next that the coefficient $q$ can not be arbitrary.
Indeed, considering $\eta>>1$ and  writing $g=1-w$ we obtain from (\ref{18}), (\ref{19}), (\ref{20})
$$
\left(\frac{d w}{d\eta}\right)^2=q\,w^2.
$$
Thus,   the function $1-g$ vanishes with an exponential rate if and only if
\be
q>0. \label{21}
\ee
The subsequent analysis depends on the parameter $\alpha$ value. We will consider separately two possibilities:
\be\label{8a}
 \alpha>0
\ee
and
\be\label{8b}
 \alpha=0.
\ee
\subsection{The case  $\alpha>0$}
Let us stress that the right-hand site in (\ref{18}) is not well defined yet since $q=q(V)$ and $V$ remains unknown.
To find $V$ we combine the second equalities in (\ref{11}) and (\ref{20}), and conclude
\be
V=\frac{1}{\alpha^2}\left\{\frac{c_3A}{1-g_*^r}-\gamma\right\}.\label{2222}
\ee
This and (\ref{12}) imply
\begin{equation}\label{221}
q=\theta-\xi(1-g_*^r),\quad \theta=c_3/\alpha^2c_1,\quad \xi=\gamma_\alpha/\alpha^2c_1A,\quad \gamma_\alpha=\gamma+\alpha^2c_0.
\end{equation}
Combining (\ref{19}) and (\ref{118}) with (\ref{221}) yields
\begin{equation}\label{222}
F\big(g_*,q(g_*)\big)\ipd \rho_1g_*^2-\rho_2g_*^{2-r}+\rho_3g_*^{2-2r}+\rho_4g_*^r-C_1=0,
\end{equation}
where
\begin{align}
&\rho_1=(1-r)(2-r+2\xi),\quad\rho_2=2(1-r)(2-\theta)+(4-3r)\xi,\label{223}\\
&\rho_3=(2-r)(1-\theta+\xi),\quad \rho_4=r\xi,\quad C_1=r(r-\theta+\xi).\notag
\end{align}
The equation (\ref{222}) has a root $g_*\in[0,1)$ if and only if  $C_1\geq0$.
However, for $g_*=0$
\begin{equation}\label{222a}
\omega'|_{\eta=+0}=-\frac{r}{p}g^{r-1}g'_\eta\big|_{\eta=+0,\,g_*=0}=-r/p\neq0,
\end{equation}
since $q(0)=r$. At the same time, the equality $F\big(g_*,q(g_*)\big)=0$ implies for each $g_*\neq0$
\be
\omega'|_{\eta=+0}=0.\label{224}
\ee
 Thus, the  assumption
\begin{equation}\label{23}
-\gamma_\alpha g_*^r<c_3A-\gamma_\alpha<rc_1A\alpha^2
\end{equation}
guaranties  the fulfilment of both the conditions  (\ref{224}) and (\ref{21}). In turn, the restriction  $g_*\in(0,1)$ implies the inequality $p<1$.
Thus we obtain the conclusion
\begin{teo}
Under the assumptions (\ref{8}), (\ref{8a}) we assume the fulfilment of the condition (\ref{23})
and define  the velocity $V$ by the formula (\ref{2222}). Then the equation (\ref{1})  has a nontrivial smooth solution (\ref{5}) with the properties (\ref{6}), (\ref{7}).
\end{teo}
\textbf{Example 1.}
 When $c_2=c_3/2$ or $c_2=c_3$ (like for the CH and  DP equations), or if $c_2=c_3/4$, or $c_2=3c_3/2$,    the function $F(g,q)$ (\ref{19}) is an algebraic   polynomial of a degree less or equal to 5. By taking into account the root $g_*=1$ of the multiplicity 2, we obtain the possibility to solve the equation  (\ref{118}) explicitly for each constant $q$ and find $g_*=G(q)$. Next we use the equality (\ref{221}) and find the root of the equation (\ref{222}) in the implicit form $g_*=G(\theta-\xi(1-g_*^r))$.

In particular, let  $r=2/3$. Then
$$
F(g,q)=(g^{2/3}-1)^2(g^{2/3}-1+3q/2).
$$
Thus
\be\label{323}
g_*^{2/3}=1-c_3A(\gamma_\alpha+rc_1\alpha^2A)^{-1}.
\ee
Substituting now $c_1=3c_3/2\alpha^2$ and $\gamma=0$ we obtain the root of (\ref{222}) for the CH equation:
$$
g_*= \big(1+c_3A/(c_0\alpha^2)\big)^{-3/2}\quad\text{if}\quad c_0>0\quad \text{and}\quad g_*=0\quad\text{if}\quad c_0=0.
$$
Respectively, the condition (\ref{23}) is satisfied for $c_0>0$ and it is broken for $c_0=0$. In the last case
$\om'_\eta|_{\eta=0}\neq 0$, therefore $\om(\eta)$ is a continuous function only.
\begin{figure}[H]
\centering
\includegraphics[width=6cm]{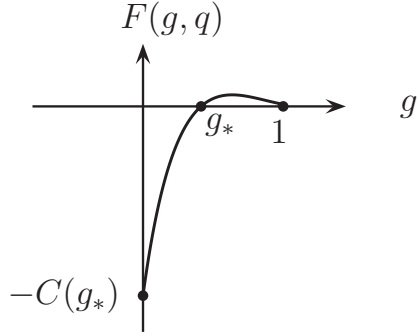}
\caption{Behavior of the function $F(g_*,q(g_*))$ for  $r=5/6$}
\label{f1}
\end{figure}
\textbf{Example 2}
Let $c_3=5c_2$, so that $r=5/6$ and $F(g_*,q(g_*))$ is a polynomial  of  degree 12. Setting $c_0=c_2=1$, $c_1=3$, $\alpha=2$, and $\gamma=0$, we solve the equation (\ref{222}) numerically using the fourth order Runge-Kutta method. For $A=1.2$ we find the desired root $g_*\approx 0.5070$ (see Fig.1).

\textbf{Example 3}
Now let $c_0=\gamma=0$ and $\alpha^2c_1>c_2+c_3$. Then $q=c_3/\alpha^2c_1<r$ and $g_*$ doesn't depend on $V$. Thus
\begin{equation}
V=c_3A/\{(1-{g_*}^r)\alpha^2\}.\label{26}
\end{equation}
\begin{remark}
Formula (\ref{222a}) shows that the soliton (\ref{5}) turns out to be a peakon in the limiting case $C_1=0$, that is for
\be
c_3A-\gamma_\alpha=rc_1A\alpha^2.\label{26b}
\ee
In accordance with Eq.(\ref{2222}), the phase velocity of this wave is
\be
V=c_0+rc_1A.\label{26c}
\ee
\end{remark}
\subsection{The case  $\alpha=0$}
Formulas (\ref{11}) and (\ref{12}) imply now the equalities
\begin{equation}
p=c_3A/\gamma, \quad q=c_3(V-c_0)/\gamma c_1.\label{26a}
\end{equation}
Thus, the supposition  (\ref{17}) allows us to define the initial value
\be
g_*=(1-p)^{1/r} \label{27b}
\ee
under the condition (\ref{16}). Consequently, instead of (\ref{222}) we obtain the following linear equation:
\be
q\mathfrak{G}=(1-r)\mathfrak{F},\label{28b}
\ee
where
\be
\mathfrak{G}=1+2p\frac{1-r}r-g_*^{2(r-1)}, \; \mathfrak{F}=\frac{r}{1-r}+2p+p^2\frac{2-r}r-\frac{r}{1-r}g_*^{2(r-1)}.\notag
\ee
\begin{lem}
Under the condition (\ref{17}) the  equation (\ref{28b})  has a solution $q_*>0$.
\end{lem}
To prove the statement it is enough to note that
\be
\mathfrak{G}|_{p=0}=0,\quad  \mathfrak{F}|_{p=0}=0,\quad \mathfrak{G},\, \mathfrak{F}\to-\infty \quad\text{as}\quad p\to1,\notag
\ee
and $\mathfrak{G}'_p<0$, $\mathfrak{F}'_p<0$ uniformly in $p\in(0,1)$.

In accordance with (\ref{26a}) we define the velocity of the soliton (\ref{5})
\be
V=c_0+q_*\gamma c_1/c_3.\label{29a}
\ee
We have thus established
\begin{teo}
Under the assumptions (\ref{8}), (\ref{8b}) we assume the fulfilment of the condition  (\ref{17})
and define  the velocity $V$ by the formula (\ref{29a}). Then the equation (\ref{1})  has a nontrivial smooth solution (\ref{5}) with the properties (\ref{6}), (\ref{7}).
\end{teo}
\begin{remark}
If
\be
A=\gamma/c_3,\label{26d}
\ee
then $p=1$ and $g_*=0$. Thus,  the phase velocity of the peakon is
\be
V=c_0+r\gamma c_1/c_3.\label{26e}
\ee
\end{remark}
Note that the equalities (\ref{26d}), (\ref{26e}) coincide with (\ref{26b}), (\ref{26c}) in the case $\alpha=0$.
\begin{remark}
Weak asymptotics constructed in \cite{Om} shows that soliton type solutions of (\ref{1}) collide elastically in the leading term with respect to $\ve<<1$ and under some additional assumptions. Results of direct numerical simulations depicted in Figures 2 and 3 confirm this conclusion. The finite-difference scheme is based on the ideas described in \cite{GO}.
\end{remark}
\begin{figure}[H]
\centering
\includegraphics[width=13cm]{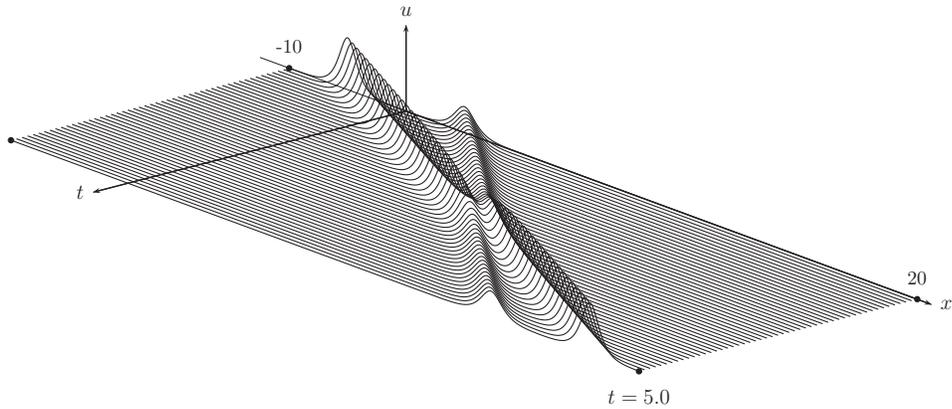}
\caption{Collision of gDP solutions in the case $\alpha=2$, $c_3=5c_2$, $c_0=c_2=1$, $c_1=3$, $\gamma=0$, $A_1=1.2$, $A_2=0.6$}
\label{f2}
\end{figure}
\begin{figure}[H]
\centering
\includegraphics[width=13cm]{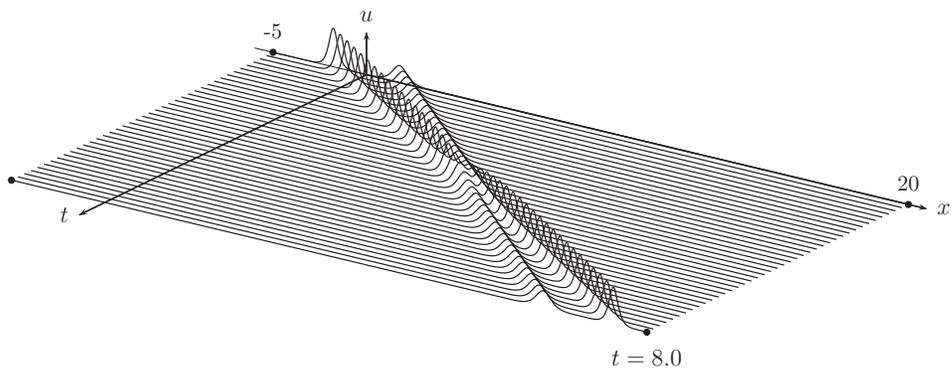}
\caption{Collision of gDP solutions in the case $\alpha=0$, $c_3=4c_2$, $c_0=c_1=c_2=1$, $\gamma=10$, $A_1=0.9$, $A_2=0.4$}
\label{f3}
\end{figure}
\section{Peakon}
Let us recall that ``peakon" is a continue  solitary wave with a jump of the first derivative. Since such functions are distributions, we write the general Degasperis-Procesi equation in the form (\ref{05}) and treat it in the weak sense. Now we define the notation
\begin{equation} \label{28}
[f]=f_+(\eta)-f_-(\eta), \quad [f]|_0=f_+(\eta)|_{\eta\to+0}-f_-(\eta)|_{\eta\to-0},
\end{equation}
and write the ansatz
\begin{equation}\label{29}
u(x,t,\ve)=A\{\om_-(\eta)+[\om]H(x-Vt)\}|_{\eta=\beta(x-Vt)/\ve},
\end{equation}
where $H$ is the Heaviside function, $H(z)=0$ for $z<0$, and $H(z)=1$ for $z>0$; the amplitude $A>0$, the phase velocity $V$, the auxiliary parameter $\beta=\sqrt{c_1/rc_3}$, and the functions $\om_\pm=\om_\pm(\eta)$ have the same sense as in (\ref{5}), but we assume now:
\begin{align}
&\om_\pm|_{\eta=\pm0}=1,\label{30}\\
&\om_+(\eta)\in C^\infty(\mathbb{R}^+), \quad \om_+(\eta)\to0\quad as\quad \eta\to+\infty,\label{31}\\
&\om_-(\eta)\in C^\infty(\mathbb{R}^-), \quad \om_-(\eta)\to0 \quad as\quad \eta\to-\infty.\label{32}
\end{align}
Obviously, (\ref{30}) implies that $[\om]|_0=0$. We assume also that the functions $\om_\pm$ are extended on $\mathbb{R^\mp}$ in a smooth manner.

Note next that $H^2=H$, thus
\begin{equation}\label{33}
u^2(x,t,\ve)=A^2\{\om_-^2(\eta)+[\om^2]H(x-Vt)\}|_{\eta=\beta(x-Vt)/\ve}.
\end{equation}
Furthermore,
\begin{align}
\ve\frac{\pa u(x,t,\ve)}{\pa x}&=A\{\om'_-(\eta)+[\om']H(x-Vt)\}|_{\eta=\beta(x-Vt)/\ve},\label{34}\\
\ve^3\frac{\pa^3 u(x,t,\ve)}{\pa x^3}&=A\{\om'''_-(\eta)+[\om''']H(x-Vt)\notag\\
&+\ve[\om'']|_0\delta(x-Vt)+\ve^2[\om']|_0\delta'(x-Vt)\}|_{\eta=\beta(x-Vt)/\ve},\notag
\end{align}
where primes denote the derivatives with respect to $\eta$, and $\delta$ is the Dirac delta-function.
Substituting (\ref{33}), (\ref{34}), and similar relation for  $u_x^2$ into (\ref{05}), and using the notation (\ref{11}), (\ref{12}), we obtain the equality
\begin{align}
\big\{\mathfrak{W}_-+[\mathfrak{W}]H\big\}&+\ve^2\big\{ [W']|_0-\frac12[(W^2)']|_0\big\}\delta' \notag\\
&+\ve\big\{ [W'']|_0-\frac{c_2-c_3}{c_3}[(W')^2]|_0-\frac12[(W^2)'']|_0\big\}\delta=0,\label{36}
\end{align}
where
\begin{equation}
\mathfrak{W}\pm=\frac{d}{d\eta}\Big\{rW_\pm^2 -rqW_\pm+W_\pm''-\frac{c_2-c_3}{c_3}(W_\pm')^2
-\frac12\big(W_\pm^2\big)''\Big\}.\label{37}
\end{equation}
Recall that the distributions $H$, $\delta$, and $\delta'$ are linearly independent. Thus by virtue of  (\ref{11}), (\ref{30}), and (\ref{36}), we deduce
\begin{equation}\label{38}
(1-p)[W']|_0=0, \quad (1-p)[W']|_0-\frac{c_2}{c_3}[(W')^2]|_0=0.
\end{equation}
Clearly, for peakons we conclude:
\be
p=1, \quad W_-(\eta)=W_+(-\eta)\quad \text{for}\quad \eta\leq0.\label{40}
\ee
Consequently, (\ref{36}) implies the problems similar to (\ref{13})
\begin{align}
&(1-W_\pm)\frac{d^2 W_\pm}{d\eta^2}=\frac{1-r}{r}\left(\frac{d W_\pm}{d\eta}\right)^2+r(q\, W_\pm-W_\pm^2),\quad \eta\in \mathbb{R}_\pm^1,\label{41}\\
&W_\pm|_{\eta=\pm0}=1, \quad W_\pm(\eta)\to0\quad as\quad \eta\to\pm\infty.\notag
\end{align}
Passing now to the equation (\ref{18}) we note that for $g_*=0$
 \be
F(g,r)=g^2(g^{-r}-1)^2.\label{118a}
\ee
Integrating we obtain the basic peakon solution:
\be
\om_\pm= \exp{(\mp r\eta)}.\label{118b}
\ee
Moreover, we obtain the same relations between $A$ and $V$ which have been described in Remarks 1 and 2. To continue note that the equalities (\ref{26b}) and (\ref{26c}) imply the conclusion: if
\be
c_3\neq rc_1\alpha^2, \quad c_0\geq0, \quad \gamma\geq0,\quad \text{and}\quad \gamma_\alpha>0,\label{118e}
\ee
then both, $V$ and $A>0$, are uniquely defined. At the same time, if
\be
c_3=rc_1\alpha^2  \quad \text{and}\quad \gamma_\alpha=0,\label{118f}
\ee
then the peakon  can be of  arbitrary amplitude.
Therefore, we conclude:
\begin{teo}
Let the conditions (\ref{8}), (\ref{26b}) be satisfied. Then the  equation (\ref{05}) has a peakon solution which propagates with the velocity  (\ref{26c}).
\end{teo}
\section{Conclusion}
Let us fix a set of structural constants $\alpha$, $c_0,\dots,c_3$, $\gamma$. Then Theorems 1 and 2 imply the existence of  the one-parametric family  of solitons  (\ref{5}), where $\om(\eta,A)$ vanishes with exponential rates and the free parameter $A$ should satisfy the restrictions
\begin{align}
 &(1-g_*^r)\frac{\gamma_\alpha}{c_3}<A<\frac{\gamma_\alpha}{c_3-r\alpha^2c_1} \quad\text{if}\quad c_3> r\alpha^2c_1\quad\text{and}\quad \gamma_\alpha>0,\label{43}\\
&(1-g_*^r)\frac{\gamma_\alpha}{c_3}<A \quad\text{if}\quad c_3= r\alpha^2c_1\quad\text{and}\quad \gamma_\alpha>0,\label{44}\\
&A=0 \quad\text{if}\quad  \gamma_\alpha=0.\label{45}
\end{align}
The soliton velocities $V$ are given by the formulas (\ref{2222}) and (\ref{29a}) for $\alpha>0$ and $\alpha=0$ respectively.
\begin{figure}[H]
\centering
\includegraphics[width=13cm]{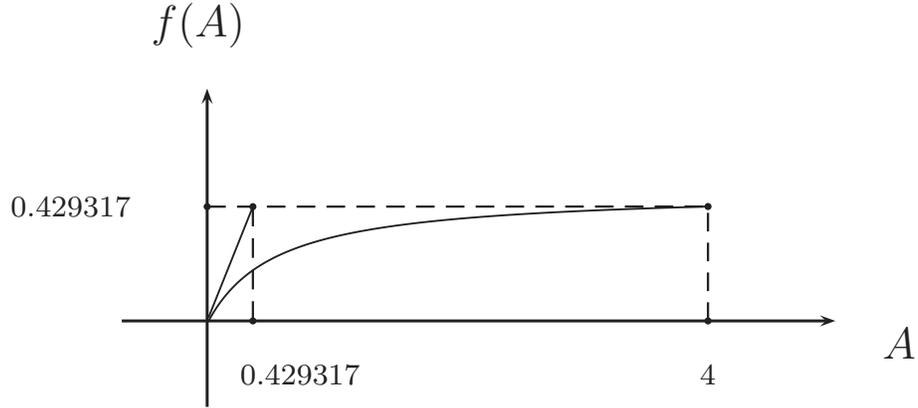}
\caption{Graphics $f=A$ and $f=\Psi(A)$ for $r=5/6$}
\label{f4}
\end{figure}
Note that  (\ref{44}) means the nonexistence of any upper bound for $A$, whereas (\ref{45}) denotes the nonexistence of solitons in this case. The lower bound in (\ref{43}), (\ref{44}) is more complicated since the root $g_*$ depends on $A$. Sometimes, there is not any actual restrictions. Indeed, for  $r=2/3$ the root $g_*$ has the form (\ref{323}) and the lower bound in (\ref{43}), (\ref{44}) is satisfied for all $A>0$. The same is true for $r=5/6$ and the parameters $\alpha,\dots,\gamma$ described in Example 2 (see Fig.4, where $\Psi(A)=\gamma_\alpha(1-g_*^r(A))/c_3$). At the same time, for $r=1/2$  is easy to prove that the restriction $\Psi(A)<A$ implies the inequality $A<4\gamma_\alpha/3c_3$. Let us note also that if   $\Psi(A)=A$, then the wave (\ref{5})  exists, but vanishes like $\sim\eta^{-2}$.

In contrast to solitons, peakons, generally speaking,  are unique waves. By virtue of Theorem 3, under the assumption (\ref{118e})  peakon should have the  amplitude
\begin{equation}
 A=\gamma_\alpha c_3/(c_3-r\alpha^2c_1), \label{47}
\end{equation}
fixed by the set of structural constants, and propagate with the fixed velocity (\ref{26c}). There is only one special case (\ref{118f}) when the general Degasperi-Procesi equation has the family of peakons with arbitrary amplitudes. Note finally that the both Degasperi-Procesi and Camassa-Holm (with $c_0=0$) equations satisfy the condition (\ref{118f}).


\begin{thebibliography}{99}

\bibitem{DegProc}
 A. Degasperis, M. Procesi,
 Asymptotic integrqability, in:  A. Degasperis, G. Gaeta (Eds.), Symmetry and Perturbation Theory, World Sientific, 23--37, 1999.

\bibitem{BBM}
 T. Benjamin, J. Bona, J. Mahony,  Model equations for long waves in nonlinear dispersive systems, Philosophical Transactions of the Royal Society of London. Series A, Mathematical and Physical Sciences, \textbf{272}, 47--78, 1972.

\bibitem{CH}
 R. Camassa, D. Holm,
An integrable shallow water equation with peaked solitons,
 Phys. Rev. Lett., \textbf{71},  1661--1664, 1993.


\bibitem{ELY}
J. Esher, Y.  Liu, Z.  Yin, Global weak solutions and blow-up structure for the Degasperis-Procesi equation,
Journal of Functional Analysis, \textbf{241}:2,  457--485, 2006.

\bibitem{BPS}
 J. Bona, W. Pritchard, L. Scott, (1980), Solitary-wave interaction, Physics of Fluids, \textbf{23}:3, 438--441, 1980.

\bibitem{DHH}
A. Degasperis, D.D. Holm, A.N.W. Hone, A new integrable equation with peakon solutions,
Theoretical and Mathematical Physics, \textbf{133}:2, 1463--1474, 2002.

\bibitem{ZQQ}
Zhijun Qiao,
M-shape peakons, dehisced solitons, cuspons and new 1-peak solitons for the Degasperis-Procesi equation, Chaos Solitons and Fractals, \textbf{37}:2, 501--507, 2008.

\bibitem{Con}
A. Constantin, On the scattering problem for the Camassa-Holm equation,
Proc. Roy. Soc. London Ser. A, \textbf{457}, 953--970, 2001.

\bibitem{LS}
H. Lundmark, J. Szmigielski, Multi-peakon solutions of the Degasperi-Procesi equation,
Invers Problems, \textbf{19}, 1241--1245, 2003.

\bibitem{Mat}
Y. Matsuno, Multisoliton solutions of the Degasperi-Procesi equation and their peakon limit,
Invers Problems, \textbf{21}, 1553--1570, 2005.

\bibitem{Om}
G. Omel'yanov, Soliton dynamics for the general Degasperis-Procesi equation,  http://arxiv.org/abs/1712.04410,  1--13,  2017.

\bibitem{GO}
M. Garcia Alvarado, G. Omel'yanov, Interaction of solitary waves for the generalized KdV equation,
Communications in Nonlinear Science and Numerical Simulation, \textbf{17}:8, 3204--3218, 2012.


\end{thebibliography}
\end{document}